\documentclass[12pt,a4paper,english]{article}
\usepackage{eurosym}

\title{\bf Zeno at the Steering Wheel}

\author{Germano D'Abramo\\
{\small Istituto Nazinale di Astrofisica, Roma, Italy}\\
{\small E--mail: {\tt Germano.Dabramo@iaps.inaf.it}}}

\date{}

\begin{document}

\maketitle

\begin{abstract}

This paper collects some reflections about an apparent incongruity 
between the usual (third-person) understanding of the probability of an 
event calculated for an extended period of time in the future (e.g.,~the 
expected probability of a driver to meet with a car accident in the next 
$M$ years) and the subjective perception of the same probability/risk 
that the person involved in that event has, instant by instant during 
that period of time. Similarities with the classical Zeno's paradoxes 
come to mind.\\

\noindent {\bf Keywords}\hspace{0.3cm} Probability $\cdot$ Time $\cdot$ 
Risk $\cdot$ Zeno's paradoxes

\end{abstract}

\section{Introduction}

In this paper we want to draw reader's attention on a question which 
could originate when one deals with the usual (third-person) treatment 
of the probability of an event calculated for an extended period of time 
in the future. More specifically, we want to focus more on the relation 
between the cumulative aspect that appears when we estimate the 
probability of an event for an extended period of time in the future and 
the subjective perception of such a probability by who is involved in 
that event and is going to experience the whole time interval in first 
person, instant by instant.

As a study case, we take into account the probability of an event which 
is in some way `close' to the experience, if not daily, of almost 
everyone of us: the probability to meet with a car accident per time 
spent at the steering wheel.

As it will be clear, all this appears to enliven the spirit of the long 
known Zeno's paradoxes.

\section{Probability of an accident}

Let~$f$ be the frequency to meet with a car accident per unit time spent 
at the steering wheel.
Actually, we do not know whether this frequency has ever been calculated 
for car driving (surely we expect it to be different from place to 
place), but it should be available for airline flights. The main problem 
with $f$ is that it is actually difficult to measure the time spent at 
the steering wheel by a single person, while it appears easier to record 
the whole time spent in flight by a passenger in his/her entire life.
  
An apparently easier option to $f$ is the frequency to meet with a car 
accident per unit kilometer traveled. In this case it is enough to add 
up the overall kilometers traveled by a sample of cars (recorded in 
their odometers) and divide the number of car accidents occurred to the 
drivers by the previous number\footnote{According to some recent 
statistics, in the safest nations of the world there is an average of 
less than one casualty per $100$ million person-kilometre~\cite{gg}.}.

For the sake of the argument, let us stick to $f$. We are not able to 
assign a numerical value to it, but it does not matter, after all.

Suppose we know of a guy, $A$, who for business reasons foresees to 
spend in his car an overall period of time equal to $T$ in the next $M$ 
years\footnote{Please note that $T$ is the sum of all the periods of 
time spent driving in $M$ years; it is not meant that the person $A$ 
drives non-stop for a time $T$.}, starting from the present time $t=0$. 
Thus, we already know that his probability to meet with {\em at least} 
one car accident in the time $T$ over $M$ years is:

\begin{equation}
P_{T}=1-e^{-fT}\qquad(\simeq fT, \quad \textrm{if}\quad fT\ll 1).
\label{eq1}
\end{equation}

The above relation is usually obtained as follows. If we consider a 
suitably short interval of time $\Delta t$ spent on the road (the sense 
of the adjective `short' will be discussed later), we can fairly assume 
that the probability $P_{\Delta t}$ to meet with a car accident in the 
interval $\Delta t$ is directly proportional to $\Delta t$ as follows,

\begin{equation}
P_{\Delta t}\simeq f\Delta t.
\label{eq2}
\end{equation}

Nothing forbids to take $\Delta t$ equal to a fraction of $T$, namely 
$\Delta t=T/N$, where $N$ is an arbitrarily large integer (obviously, it 
holds $N\cdot\Delta t=T$).

Now, since the probability of {\bf not} having a car accident in the 
time interval $\Delta t$ is $1-f\Delta t$, the probability of {\bf not} 
having a car accident in the period $T=N\cdot\Delta t$ is equal to,

\begin{equation}
\overline{P}_{T}=(\overline{P}_{\Delta t})^N=(1-f\Delta t)^N=
\biggl(1-\frac{fT}{N}\biggr)^N.
\label{eq3}
\end{equation}

In writing eq.~(\ref{eq3}), a further assumption has been made: the 
probability of {\bf not} having a car accident in the {\em nth} time 
interval $\Delta t$ is independent of the (non) occurrence of the same 
event in the previous overall time $(n-1)\cdot\Delta t$, obviously 
involving the same driver and with $1\leq n \leq N$, much like the event 
of not getting a specific number in rolling a die is independent of not 
getting the same number in previous throws of the same die (stochastic 
independence). All this allows us to write $\overline{P}_{T}$ as the 
product of all the $(1-f\Delta t)$ terms.

Now, as $N$ approaches infinity, equation~(\ref{eq3}) approaches 
$e^{-fT}$, and since $P_{T}=1-\overline{P}_{T}$, we have exactly that 
$P_{T}=1-e^{-fT}$.

By using Maclaurin series, we can verify how~$1-e^{-fT}$ can be well 
approximated by $fT$, if $fT\ll 1$,

\begin{equation}
1-e^{-fT}=fT-\frac{1}{2}(fT)^2+\frac{1}{6}(fT)^3-\cdots\approx fT,
\qquad \textrm{if $fT\ll 1$}.
\label{eq4}
\end{equation}

It is clear from eq.~(\ref{eq1}) and~(\ref{eq4}) that the longer is $T$, 
the higher is the probability for a person to be involved in a car crash.

What is written so far is perfectly understandable and straightforward. 
It is the way in which such things are handled in actuarial practice (if 
one travels more kilometers per year or, equivalently, spends more time 
at the steering wheel per year, then he is considered more at risk and 
maybe he must pay more). If we know in advance that person $A$ will 
spend driving an overall time $T$, then we already know that his/her 
future (overall) risk to be in a crash is larger than the risk taken by 
a person $B$ that, for instance, will spend an overall shorter time 
$T'<T$ at the steering wheel.

The previous approach is perfectly meaningful from a perspective which 
is, let's say, a `third-person' one. But what can be said from a 
`subjective' point of view? Namely, from the point of view of the person 
that is going to take the risk. From a subjective point of view things 
seem to behave differently.

In the aforementioned scheme, the past time spent in a car (with or 
without accidents) does not affect the present or future risk to be in a 
car crash, much like the outcome of ten heads in a row in a coin-tossing 
does not affect the probability that the eleventh outcome will or will 
not be a tail (it is again the stochastic independence used to derive 
eqs.~(\ref{eq1}) and (\ref{eq4})).

Thus, whenever $A$ has spent some of his/her time in the car without 
accidents and is about to spend some other time, the subjective 
perception is that the risk restarts each time from the present; in that 
very same instant, every instant after another, everything starts from 
scratch and the past does not count.

To be clear, when we use the term `instant' we do not refer to an ideal, 
dimensionless {\em time point}, reminiscent of a point on the real 
number line. Instead, we mean to refer in any case to a {\em time 
interval}. When we talk about the present instant, we mean a time 
interval, arbitrarily small, that an individual subjectively and 
consciously feels as {\em his/her own present}. In the next Section, 
this issue will be discussed more thoroughly.

At this point, some quite natural questions come to mind. Why should $A$ 
feel in (his/her) every present instant that he/she is overall (namely, 
if one takes into account the whole driving time $T$ foreseen in the 
coming $M$ years) more at risk than who spends less time on the road? If 
the frequency of $A$ to be in a car crash is $f$ and if, instant by 
instant, $f$ depends neither upon the (peculiar) past time nor upon the 
future that is to come\footnote{By definition, we must consider $f$, 
once derived e.g.~through statistical averaging, as a constant number.}, 
why does the overall risk\footnote{In the rest of the paper we use the 
words `risk' and `probability' indistinguishably. Strictly speaking, 
they are not, `risk' being the probability of an event times the 
resulting loss or cost. For the sake of simplicity, here we take as loss 
or cost parameter an adimentional number between 0 and 1 and assume that 
it is always 1. This way `risk' and `probability' are the same.} of $A$ 
{\em add up} to $1-e^{-fT}$?

This appears to have the flavor of a straight Zenonian paradox applied 
to probability.

To sum up, the frequency per unit time to meet with a car accident is 
$f$. The past time already spent driving does not affect the risk to be 
in a car crash for the present. Then, why does the overall risk (the 
risk for the overall foreseen driving time $T$ in the coming $M$ years 
and estimated `from outside' at the beginning of the period $[0;M]$, 
$t=0$ being the beginning of the period under analysis) add up to 
$1-e^{-fT}$, giving to $A$, at the beginning of the period $[0;M]$, the 
uneasy feeling that he/she is going to take an overall greater risk?

As a matter of fact, the risk assessment described above, 
eqs.~(\ref{eq1}) and (\ref{eq4}), seems to have a meaning and to provide 
a consistent picture only if one consider the driving experience of 
person $A$ from outside and for the overall time that $A$ is going to 
spend on the road. But if one takes into account the subjective 
perspective, that of the person involved in the risk and in (his/her) 
every present instant, that assessment seems to have no longer sense.

\section{Further remarks on the subjective perception of probability}

In order to clarify the above point once more, let us focus a bit on the 
concept and definition, far from trivial, of {\em time instant}.

The common perception that everyone of us has of the passage of time is 
that it is a {\em seamless, continuous flow}. Thus, we think about time 
as being infinitely divisible, much like mathematicians describe the 
real number line. Such a view instinctively belongs to the human being, 
it is inborn. Moreover, every conscious meditation on the topic seems to 
reinforce this belief.

Nevertheless, many ancient thinkers, and many modern scholar as well 
(with renewed vigor after the development of quantum mechanics), 
suggested the possibility that time is not infinitely divisible and that 
its flow is discretized.

As a matter of fact, it is definitely not easy to image the existence of 
a smallest interval of time, below which it is not possible to 
physically conceive a flow of time. However, if time turns out to be not 
an intrinsic and independent entity, but a relative physical quantity, 
namely definable only through close relation to the physical processes 
that occur in physical reality, then quantum mechanics seems to suggest 
that below a time unit, the Planck time $t_P\simeq 
10^{-43}\,\textrm{sec}$, time does not exist and it makes no sense to 
talk about what happens within $10^{-43}\,\textrm{sec}$. 
Oversimplifying, time may exist only because there are physical systems 
that through repetitions of one or another standard cyclical event allow 
to define and measure it --operational definition-- and it makes no 
sense to talk about time beyond the limits of the physical instruments 
used to measure these repetitions --e.g.~Heisenberg uncertainty 
principle.

Here we stick to the conception that time is a seamless, continuous flow 
and that it is infinitely divisible. What we are putting forward in the 
paper, however, still holds even with the existence of a smallest time 
unit (e.g., $t_P$).

Then, suppose that the {\em present instant} can be defined as an 
infinitesimal time interval, $dt$, felt as belonging to our present.
According to the mathematical definition, the infinitesimal is a number 
with an absolute value greater than zero, yet it is less than any 
positive real number. A number~$x>0$ is infinitesimal if and only if 
every sum $x + ... + x$ of a finite number of terms is always less than 
any positive number, no matter how big is the finite number of terms.

Thus, going back to the main topic, the probability to be in a car crash 
that everyone of us {\em feels} in every {\em present instant} is:

\begin{equation}
P_{\cal P}=fdt, 
\label{eq5}
\end{equation}
namely, it is an infinitesimal quantity (since $f$ is a finite 
quantity). Therefore, in every present instant it is practically zero, 
since it is always less than any finite real number. It would seem, 
therefore, that nothing can happen to us.

Moreover, $P_{\cal P}$ does not depend on,

\begin{equation} 
\int_0^{T_{\cal P}} fdt, 
\label{eq6} 
\end{equation}
namely, according to what we have said before, it does not depend upon 
the past (where $T_{\cal P}$ is the present time, the upper bound of the 
overall period $[0;T_{\cal P}]$ spent at the steering wheel so far).

Therefore, from a subjective perspective not only the probability to 
meet with a car accident in every present instant is always $fdt$ (it 
does not dependent on the past), but it always has an infinitesimal 
value.

This result is, prima facie, one that our minds do not want to accept, 
much like Zeno's arguments on motion (see, for instance, the Dichotomy 
Paradox \cite{boy,sep}).

About the independence of $P_{\cal P}$ from past events, let us give 
some real life examples that should help understand the point.

Consider the case of collection of money for a day's work: if I get paid 
$X$~{\euro} an hour for my work, the total amount of money I get in a 
day's work (e.g.~$8$ hours) won't be $X$, but obviously $8\times 
X$~{\euro}. Money is something material that accumulates and its amount 
depends obviously on the elapsed time and on the future time over which 
one expects to receive it (e.g., at the end of the week I already know 
that I will be paid $5\times 8\times X$~{\euro}). Conversely, as has 
been shown before, the probability $fdt$ does not depend upon both the 
past and the future and it is always the same in every present instant.

Within $m$ hours I will have an amount of money equal to $\approx 
m\times X$~{\euro}, while in the next $m$ hours spent on the road I feel 
a probability to be in a car accident, in every `present instant' of 
these $m$ hours, always equal to $fdt$, although the probability to meet 
with a car accident in the {\em overall} next $m$ hours, if it is 
calculated before I start driving and looking at the near future, is 
$1-e^{-m\times 3600\,\textrm{\scriptsize sec}\times f}$ (if $f$ is 
measured in sec$^{-1}$). It can be fairly approximated to $\approx 
m\times 3600\,\textrm{sec}\times f$, since usually $m\times 
3600\,\textrm{sec}\times f\ll 1$.

An example in which this sort of `cumulative' aspect of the probability 
of an event calculated for an extended period of time in the future has 
a plain meaning is the probability to get sick due to exposure to 
dangerous and poisonous agents.

In this case, the time already spent in contact with a dangerous 
chemical substance or under the exposure of e.m.~waves is important, if 
not crucial, in the assessment of the future risk. It is perfectly 
conceivable that the exposure to the agents cumulatively affects the 
overall risk, since chemical agents or ionizing radiations have effects 
that accumulate at cellular level and traces of past exposure remain in 
the biological tissue (impairment, genetic mutation, damage). In the 
case of the probability of a car crash, or in every case of probability 
of an event independent of past occurrences, nothing accumulates. 
Nevertheless, in taking into account the overall risk from a 
third-person perspective and at the beginning of the whole period under 
study, see eq.~(\ref{eq1}) and (\ref{eq4}), it seems like it accumulates 
with time.

Again, the overall risk of a car accident for the coming $m$ hours, 
reckoned now and looking at the future, is $\approx m\times 
3600\,\textrm{sec}\times f$, but the risk felt in every present instant 
(infinitesimal $dt$) at the steering wheel {\em during those hours} is 
always $fdt$ and it is infinitesimal, namely practically zero.

\section{Conclusions}

Like it happens with Zeno's Paradoxes, the argument shown in the paper 
appears prima facie convincing and even sound, but we do know that it 
cannot be true. This is the core of any logical paradox. Perhaps, it may 
suggest an incongruity between reality and our mathematical treatment of 
it or, maybe, a conflict between objective reality (and its mathematical 
treatment) and our instinctive and subjective perception of it. In any 
case, our understanding of the flow of time is also called into 
question.

With Zeno we would say that like an arrow thrown toward a target can 
never reach it, so every event having a non-zero frequency per unit time 
to happen will actually not happen to us. In the physical reality, 
however, arrows reach quickly their target (if thrown properly) and 
events, good and bad, happen around and to us all.

\end{document}